\documentclass{article}
\usepackage{amsfonts}
\usepackage{amsmath}
\usepackage{enumerate}
\usepackage{amsthm}
\usepackage{hyperref}

\begin{document}

\newtheorem{Def}{Definition}[section]
\newtheorem{Thm}[Def]{Theorem}
\newtheorem{Prop}[Def]{Propostion}
\newtheorem{Rk}[Def]{Remark}
\newtheorem{Lma}[Def]{Lemma}
\newtheorem{Cor}[Def]{Corollary}

\numberwithin{equation}{section}

\title{Generalized maximum principles and stochastic completeness for pseudo-Hermitian manifolds}
\author{Yuxin Dong\footnote{Supported by NSFC grant No. 11771087, and LMNS, Fudan}, \ Weike Yu\footnote{Corresponding author.}}
\date{}
\maketitle
\begin{abstract}
In this paper, we establish a generalized maximum principle for pseudo-Hermitian manifolds. As corollaries, Omori-Yau type maximum principles for pseudo-Hermitian manifolds are deduced. Moreover, we prove that the stochastic completeness for the heat semigroup generated by the sub-Laplacian is equivalent to the validity of a weak form of the generalized maximum principles. Finally, we give some applications of these generalized maximum principles. 
\end{abstract}

\section{Introduction}
In 1967, Omori \cite{[Om]} first introduced a generalized maximum principle on a complete Riemannian manifold with sectional curvature bounded below, that is, on such a Riemannian manifold $N$, for any $\epsilon>0$ and $C^2$ function $u:N\rightarrow \mathbb{R}$ which is bounded from above, there exists a sequence $\{x_k\}\subset N$ such that 
\begin{align}
\lim_{k\rightarrow\infty}u(x_k)=\sup_M{u},\ \lim_{k\rightarrow\infty} |\nabla u(x_k)|=0, \  \limsup_{k\rightarrow\infty} \text{Hess}(u)(X_k,X_k)\leq 0\label{1.1}
\end{align}  
for any $X_k\in T_{x_k}N$ with $|X_k|=1$. Later, Yau generalized it to a complete Riemannian manifold with Ricci curvature bounded below (cf. \cite{[CY]},  \cite{[Ya]}), where the conclusion  \eqref{1.1} was modified to
\begin{align}
\lim_{k\rightarrow\infty}u(x_k)=\sup_M{u},\ \lim_{k\rightarrow\infty} |\nabla u(x_k)|=0, \ \limsup_{k\rightarrow\infty} \Delta u(x_k)\leq 0.
\end{align}
Since then, these maximum principles have been extended by several authors (cf. \cite{[Bo]}, \cite{[CX]}, \cite{[HS]}, \cite{[PRS-2]}, \cite{[RRS]}, \cite{[Ta]}, etc.) and have become an important analytical tool in differential geometry. Furthermore, S. Pigola, M. Rigoli and A.G. Setti \cite{[PRS-1]} proved the equivalence between the stochastic completeness of a Riemannain manifold and the validity of a weak form of the above maximum principles.

The present paper is mainly devoted to the study of generalized maximum principles and stochastic completeness for pseudo-Hermitian manifolds. For the convenience of readers, an introduction to pseudo-Hermitian geometry will be presented in Section \ref{section2}, and the definition of stochastic completeness for pseudo-Hermitian manifolds will be given in the front part of Section \ref{section4}. Roughly speaking, pseudo-Hermitian manifolds are CR manifolds of hypersurface type with positive definite pseudo-Hermitian structures (see Section \ref{section2} for the precise definition). Let $(M^{2m+1},H,J,\theta)$ denote a pseudo-Hermitian manifold, where $(H,J)$ is a CR structure of type $(m,1)$ and $\theta$ is a pseudo-Hermitian structure. It carries a positive definite Levi form $L_\theta$ on $H$ which is induced by the complex structure $J$ and the pseudo-Hermitian structure $\theta$. Furthermore, on a pseudo-Hermitian manifold, there is a natural 2-step sub-Riemannian structure $(H, L_\theta)$ which induces the Carnot-Carath\'eodory distance $r_{cc}$ on $M$. So we can see that the pseudo-Hermitian manifolds carry rich geometric structures.

In this paper, Firstly, by using a similar method as in \cite{[PRS-2]}, we establish a generalized maximum principle for pseudo-Hermitian manifolds (see Theorem \ref{theorem3.1}). Consequently, in terms of the sub-Laplacian and Hessian comparison theorems respectively, we obtain the following generalized maximum principles of Yau and Omori types respectively.
\begin{Thm}
Let $(M^{2m+1},H,J,\theta)$ be a complete pseudo-Hermitian manifold with pseudo-Hermitian Ricci curvature bounded below and $\|A\|_{C^1}$ bounded above where $A$ is the pseudo-Hermitian torsion. Then for every $C^2$ function $u:M\rightarrow \mathbb{R}$ which is bounded from above, there exists a sequence $\{x_k\}\subset M$ such that
\begin{align}
\lim_{k\rightarrow\infty}u(x_k)=\sup_M{u},\ \lim_{k\rightarrow\infty}|\nabla{u(x_k)}|=0, \  \limsup_{k\rightarrow\infty} {\Delta_b u(x_k)\leq 0},
\end{align}
where the norm $|\cdot|$ is induced by the Webster metric $g_\theta$, and $\Delta_b$ is the sub-Laplacian (see \eqref{2.18..}).
\end{Thm}
\begin{Thm}
Let $(M^{2m+1},H,J,\theta)$ be a complete pseudo-Hermitian manifold with horizontal sectional curvature $K_{\sec }^{H}$ bounded below and $\|A\|_{C^1}$ bounded above where $A$ is the pseudo-Hermitian torsion. Then for every $C^2$ function $u:M\rightarrow \mathbb{R}$ which is bounded above, there exists a sequence $\{x_k\}\subset M$ such that
\begin{gather}
\lim_{k\rightarrow\infty}u(x_k)=\sup_M{u},\
\lim_{k\rightarrow\infty} |\nabla{u(x_k)}|=0,\
\limsup_{k\rightarrow\infty} \text{Re} \nabla du(X_k, \overline{X_k})\leq 0
\end{gather}
for every $X_k\in T_{1,0}M_{x_k}$ with $|X_k|=1$, where $\text{Re} \nabla du(X_k, \overline{X_k})$ is defined as in \eqref{2.19..}.
\end{Thm}

Moreover, we study the relationship between the generalized maximum principles and stochastic completeness for the heat semigroup generated
by the sub-Laplacian. Indeed, we prove that the stochastic completeness of a pseudo-Hermitian manifold is equivalent to the validity of a weak form of the generalized maximum principles, which is a pseudo-Hermitian form of Theorem 1.1 in \cite{[PRS-1]}.
\begin{Thm}
Let $(M^{2m+1},H,J,\theta)$ be a non-compact pseudo-Hermitian manifold. The following statements are equivalent:
\begin{enumerate}[(1)]
\item $M$ is stochastically complete.
\item For every $\lambda>0$, the only nonnegative bounded $C^2$ function on $M$ with $\Delta_b  u=\lambda u$ is $u\equiv0$.
\item If $u$ is a $C^2$ function on $M$ with $\sup_M u<+\infty$, then $\inf_{\Omega_\alpha} \Delta_b u\leq 0$ for any $\alpha>0$, where $\Omega_\alpha=\{x\in M:\ u(x)>\sup_M u-\alpha \}$.
\item For any $u\in C^2(M)$ with $\sup_M u<+\infty$, there exists $\{x_n\}\subset M$ such that $u(x_n)\geq \sup_M u-\frac{1}{n}$ and $\Delta_b u(x_n)\leq \frac{1}{n}$.
\end{enumerate}
\end{Thm}

At last, we will give two applications of the generalized maximum principles for pseudo-Hermitian manifolds, which are similar to some results in Riemannian case (cf. \cite{[RRS]}, \cite{[Su]}, ect.). One is to study differential inequalities on a pseudo-Hermitian manifold, and the other is to investigate the obstruction of pseudo-Hermitian scalar curvature to CR conformal deformations.

\section{Preliminaries}\label{section2}
In this section, we present some facts and notations in pseudo-Hermitian geometry (cf. \cite{[DT]}, \cite{[BG]}).

A CR manifold is a real smooth orientable $2m+1$ dimensional differentiable manifold $M$ equipped with a complex subbundle $T_{1,0}M$ of complex rank $m$ of the complexified tangent bundle $TM\otimes \mathbb{C}$ such that
\begin{align}
T_{1,0}M\cap T_{0,1}M=\{0\}, \ \ [\Gamma(T_{1,0}M),\Gamma(T_{1,0}M)]\subseteq\Gamma(T_{1,0}M)\label{2.1}
\end{align}
where $T_{0,1}M=\overline{T_{1,0}M}$. The subbundle $T_{1,0}M$ is called a CR structure on $M$. The CR structure corresponds to a real $2m$-distribution, called Levi distribution, which is given by
\begin{align}
H=Re\{T_{1,0}M\oplus T_{0,1}M\}.
\end{align}
Clearly, $H$ carries a natural complex structure $J$, which is defined as $J(X+\overline{X})=\sqrt{-1}(X-\overline{X})$ for any $X\in T_{1,0}M$. 

Since both $M$ and $H$ are orientable, there is a global nowhere vanishing 1-form $\theta$ such that $H=\ker \theta$. Such a section $\theta$ is called a pseudo-Hermitian structure of $M$, and the corresponding Levi form is defined as
\begin{align}
L_\theta(X,Y)=d\theta(X,JY)
\end{align}
for any $X,Y\in H$. The second condition in \eqref{2.1} implies that $L_\theta$ is $J$-invariant, and thus symmetric. If the Levi form $L_\theta$ is positive definite on $H$, $(M^{2m+1}, H, J)$ is said to be strictly pseudoconvex. Such a quadruple $(M^{2m+1},H,J,\theta)$ is called a pseudo-Hermitian manifold.

The non-degeneracy of $L_{\theta }$ on $H$ implies that all sections of $H$ together with their Lie
brackets span $T_{x}M$ at each point $x\in M$. In fact, $(M^{2m+1},H,L_{\theta })$ is a $2$-step sub-Riemannian manifold. A Lipschitz curve $\gamma: [0,l]\rightarrow M$ is said to be horizontal if $\gamma'(t)\in H_{\gamma(t)}$ a.e. in $[0,l]$. From the well-known theorem of Chow-Rashevsky (cf. \cite{[Ch]}, \cite{[Ra]}), it follows that for any two points $p,q\in M$, there is a horizontal Lipschitz curve joining $p$ and $q$. Consequently, the Carnot-Carath\'eodory distance is defined by
\begin{align*}
r_{cc}(p, q)=\inf\{\int_0^l \sqrt{L_\theta(\gamma', \gamma')}\, dt\ |\ &\gamma: [0, l]\rightarrow M \ \text{is a horizontal Lipschitz curve,}\\
&\ \gamma(0)=p,\  \gamma(l)=q  \},
\end{align*}
which induces to a metric space structure on $(M^{2m+1},H,L_\theta)$.

For a pseudo-Hermitian manifold $(M^{2m+1},H,J,\theta)$, there is a unique globally defined nowhere zero tangent vector field $\xi$ on $M$, which is called the Reeb vector field, such that 
\begin{align}
\theta(\xi)=1,\ \  d\theta(\xi, \cdot)=0.
\end{align}
Consequently there is a splitting of the tangent bundle $TM$
\begin{align}
TM=H\oplus L,
\end{align}
where $L$ is the trivial line bundle generated by $\xi$. This direct sum decomposition can induce two natural projections denoted by $\pi_H:TM\rightarrow H$ and $\pi_L:TM\rightarrow L$. 

Let us extend $J$ to a $(1,1)$-tensor field on $M$ by requiring
\begin{align}
J\xi=0.
\end{align}
Set 
\begin{align}
G_\theta(X,Y)=L_\theta(\pi_HX, \pi_HY)
\end{align}
for any $X,Y\in TM$, then the $J$-invariance of $L_\theta$ implies that $G_\theta$ is also $J$-invariant. Since $G_\theta$ coincides with $L_\theta$ on $H\times H$, $L_\theta$ can be extended to a Riemannian metric $g_\theta$ on $M$ given by
\begin{align}
g_\theta=G_\theta+\theta\otimes\theta
\end{align}
which is usually called the Webster metric. 

For a pseudo-Hermitian manifold $(M^{2m+1},H,J,\theta)$, one can define two notions of completeness by using the Carnot-Carath\'{e}odory distance $r_{cc}$ of $L_{\theta }$ and the Riemannian distance $r$ of $g_{\theta}$ respectively. Actually, they are equivalent, because $r_{cc}$ and $r$ are locally controlled by each other (cf. e.g., \cite{[NSW]}).

On a pseudo-Hermitian manifold, there is a canonical connection preserving the CR structure and the Webster metric, which is called the Tanaka-Webster connection.
\begin{Thm}[cf. \cite{[Tan]}, \cite{[We]}]\label{theorem 2.1}
Let $(M^{2m+1},H,J,\theta)$ be a pseudo-Hermitian manifold with the Reeb vector field $\xi$ and Webster metric $g_\theta$. Then there is a unique linear connection $\nabla$ on $M$ satisfying the following axioms:
\begin{enumerate}[(1)]
\item $\nabla_X\Gamma(H)\subset \Gamma(H)$ for any $X\in \Gamma(TM)$.
\item $\nabla J=0$, $\nabla g_\theta=0$.
\item The torsion $T_\nabla$ of $\nabla$ is pure, that is, 
\begin{align*}
T_\nabla(X,Y)=2d\theta(X,Y)\xi\ and\  T_\nabla(\xi,JX)+JT_\nabla(\xi,X)=0
\end{align*}
for any $X,Y\in H$.
\end{enumerate}
\end{Thm}
The pseudo-Hermitian torsion, denoted by $\tau$, is a $(1, 1)$-tensor on $M$ defined by
\begin{align}
\tau(X)=T_\nabla(\xi,X)
\end{align}
for any $X\in TM$. Set
\begin{align}
A(X,Y)=g_\theta(\tau(X),Y)
\end{align}
for any $X,Y\in TM$. A pseudo-Hermitian manifold is called a Sasakian manifold if $\tau\equiv0$ (or equivalently, $A\equiv0$). Note that the properties of $\nabla$ in Theorem \ref{theorem 2.1} imply that
$\tau(T_{1,0}M)\subset T_{0,1}M$ and $A$ is a trace-free symmetric tensor field. In addition, one can define the $C^1$ norm of $A$ on a open subset $\Omega$ of $M$ as
\begin{align}
\|A\|_{C^1(\Omega)}=\max_{x\in \Omega}\{\|A\|(x),\|\nabla A\|(x)\}.
\end{align}
If $\Omega=M$, the above quantity is denoted by $\|A\|_{C^1}$ for simplicity.

The curvature tensor $R$ associated with the Tanaka-Webster connection is given by 
\begin{align}
R(X,Y)Z=\nabla_X\nabla_Y Z-\nabla_Y\nabla_X Z-\nabla_{[X,Y]}Z,
\end{align}
and set $R(X,Y,Z,W)=g_\theta(R(Z,W)X, Y)$ where $X,Y,Z,W\in TM$. Let $\{\eta_i\}_{i=1}^m$ be a local unitary frame of $T_{1,0}M$ with respect to the Hermitian structure determined by the Webster metric and let $\eta_0=\xi$. According to \cite{[We]}, $R_{ABCD}$, which are the components of $R$ under the local frame $\{\eta_0, \eta_i, \eta_{\bar{i}}\}$, have the following properties:
\begin{align}
R_{i\bar{j}k\bar{l}}=-R_{\bar{j}ik\bar{l}}=-R_{i\bar{j}\bar{l}k},\ \ \ R_{i\bar{j}k\bar{l}}=R_{k\bar{j}i\bar{l}}=R_{i\bar{l}k\bar{j}}.
\end{align}
Similar to the Riemannian case, E.Barlerra investigated the following sectional curvature of Tanaka-Webster connection in \cite{[Ba]}:
\begin{align}
K_{sec}(\sigma)=-R(u,v,u,v)\label{2.13.}
\end{align}
for any 2-plane $\sigma\subset T_xM$, where $\{u,v\}$ is a $g_\theta$-orthonomal basis of $\sigma$. In particular, if $\sigma\subset H$, the above quantity is called the horizontal sectional curvature and denoted by $K^H_{sec}(\sigma)$.
The Ricci tensor of the Tanaka-Webster connection is defined as 
\begin{align}
Ric(Y,Z)=trace\{X\rightarrow R(X, Z)Y \}
\end{align}
for any $Y,Z\in TM$. We also have the pseudo-Hermitian Ricci tensor defined by 
\begin{align}
R_{i\bar{j}}=Ric(\eta_i, \eta_{\bar{j}})=\sum_{k=1}^m R_{k\bar{k}i\bar{j}}.
\end{align}
The pseudo-Hermitian scalar curvature is given by
\begin{align}
R=\frac{1}{2}trace(Ric)=\sum_{k=1}^m R_{k\bar{k}}.
\end{align}

For a smooth function $u: (M, H, J, \theta)\rightarrow \mathbb{R}$, one can define the Hessian of $u$ with respect to the Tanaka-Webster connection $\nabla$ by
\begin{align}
(\nabla du)(Y,Z)=Y(Zu)-(\nabla_YZ)u\label{2.18...}
\end{align}
for $Y, Z\in TM$. Set $u_{k\bar{l}}=(\nabla du) (\eta_k, \eta_{\bar{l}})$, then

\begin{align}
Re \nabla du (X,\overline{X})=\frac{1}{2}\sum_{k,l=1}^m\left(u_{k\bar{l}}X^kX^{\bar{l}}+\overline{u_{k\bar{l}}X^kX^{\bar{l}}}\right),\label{2.19..}
\end{align}
where $X=\sum_{k=1}^mX^k\eta_k\in T_{1,0}M$, $X^{\bar{l}}=\overline{X^l}$. Let $\nabla^\theta$ be the Levi-Civita connection of the Riemannian manifold $(M, g_\theta)$, and the Hessian of $u$ with respect to $\nabla^\theta$ can be similarly defined as in \eqref{2.18...}. According to the relationship between $\nabla$ and $\nabla^\theta$ (cf. Lemma 1.3 in \cite{[DT]}), it is easy to show that (cf. the proof of Lemma 4.6 in \cite{[DRY]})
\begin{align}
Re \nabla du (X,\overline{X})=\frac{1}{2}\left( \nabla^\theta du (X,\overline{X})+\nabla^\theta du (\overline{X},X)\right)
\end{align}
for any $X=\sum_{k=1}^mX^k\eta_k\in T_{1,0}M$. 

Analogous to the Laplace operator on a Riemannian manifold, there is a degenerate elliptic operator on a pseudo-Hermitian manifold, which is called the sub-Laplace operator. Precisely, the sub-Laplacian of $u$ is defined by
\begin{align}
\Delta_b u=trace_H(\nabla du)=\sum_k (u_{k\bar{k}}+u_{\bar{k}k}).\label{2.18..}
\end{align}

If we choose a local $G_\theta$-orthonormal real frame $\{X_\alpha\}_{\alpha=1}^{2m}$ of $H$ defined on the open set $U\subset M$, then the sub-Laplacian of $u$ can be expressed by
\begin{align}
\Delta_b u=trace_H(\nabla du)= \sum_{\alpha=1}^{2m} X_{\alpha}^2u+X_0 u
\end{align} 
where $X_0=-\sum_{\alpha=1}^{2m}\nabla_{X_{\alpha}}X_{\alpha}\in H$. From the positive definiteness of the Levi form $L_\theta$ on $H$, it follows that $X_0, X_1,\dots,X_{2m}$ together with their commutators span the tangent spaces at any point of $U$, and therefore, $\Delta_b$ is hypoelliptic according to Theorem 1.1 in \cite{[Ho]}. Furthermore, Proposition 2.3 in \cite{[DT]} states that
\begin{align}
\int_M (\Delta_bu) v\ d\mu=-\int_M L_\theta(\nabla^Hu, \nabla^Hv)\ d\mu\label{2.20..}
\end{align}
holds for any smooth function $u,v:M\rightarrow \mathbb{R}$, at least one of compact support, where $\mu$ is the Borel measure given by the volume form $\theta\wedge(d\theta)^m$ of $M$.

At the end of this section, we recall the Folland-Stein spaces briefly (cf. \cite{[FS]}). For any $k\in\mathbb{N}_+$ and any $p$ with $1<p<+\infty$, the Sobolev space compatible to CR structure is defined by
\begin{align}
S_k^p(U)=\left\{f\in L^p(U,\mu): X_{i_1}X_{i_2}\dots X_{i_s}f\in L^p(U,\mu), s\leq k, X_{i_j}\in \{X_\alpha\}_{\alpha=1}^{2m}\right\}
\end{align}
with the norm
\begin{align}
\|f\|_{S_k^p(U)}=\|f\|_{L^p(U,\mu)}+\sum\|X_{i_1}X_{i_2}\dots X_{i_s}f\|_{L^p(U,\mu)},
\end{align}
where the sum is taken over all ordered monomials $X_{i_1}X_{i_2}\dots X_{i_s}$, $1\leq s\leq k$ of the local $G_\theta$-orthonormal real frame $\{X_\alpha\}_{\alpha=1}^{2m}$ of $H$. Such Sobolev spaces are usually called the Folland-Stein spaces. By partition of unity, we can also define $S^p_k(\Omega)$, where $\Omega$ is any open subset of $M$. Under these generalized Sobolev spaces, we can investigate the existence and regularity of the solutions of subelliptic equations. 

\section{Generalized maximum principles for pseudo-Hermitian manifolds}\label{section3}
In this section, we will establish generalized maximum principles for pseudo-Hermitian manifolds.
\begin{Thm}\label{theorem3.1}
Let $(M^{2m+1},H,J,\theta)$ be a pseudo-Hermitian manifold. Assume that there exists a non-negative $C^2$ function $\gamma$ on $M$ satisfying 
\begin{gather}
\lim_{x\rightarrow \infty} \gamma(x)=+\infty,\label{3.1}\\
|\nabla^H \gamma|\leq A\gamma^{\frac{1}{2}} \ \ \text{outside some compact subset of}\ M, \label{3.2}\\
\Delta_b \gamma \leq B\gamma^{\frac{1}{2}}G(\gamma^{\frac{1}{2}})^{\frac{1}{2}}\ \ \text{outside some compact subset of}\ M,\label{3.3}
\end{gather}
where \eqref{3.1} means that $\forall \eta>0, \exists\ \text{compact set}\ K_\eta \subset M\ \text{such that}\ \gamma(x) > \eta \ \text{whenever}\ x\notin K_\eta$, the norm $|\cdot|$ is induced by the Webster metric, $\nabla^H=\pi_H\circ\nabla$ is the horizontal gradient, and $A, B$ are two positive constants, $G$ is a smooth function on $[0, +\infty)$ with
\begin{gather}
\begin{array}{llll}
i)\ G(0)>0, & ii)\ G'\geq 0, \\
iii)\ G^{-\frac{1}{2}}\notin L^1(+\infty),  & iv)\ \limsup_{t\rightarrow +\infty} \frac{tG(t^{\frac{1}{2}})}{G(t)}<+\infty.
\end{array}
\label{3.4}
\end{gather}
Then for every $C^2$ function $u:M\rightarrow \mathbb{R}$ which is bounded from above, there exists a sequence $\{x_k\}\subset M$ such that
\begin{align}
\lim_{k\rightarrow\infty}u(x_k)=\sup_M{u},\ \lim_{k\rightarrow\infty} |\nabla^H{u(x_k)}|=0, \  \limsup_{k\rightarrow\infty} {\Delta_b u(x_k)\leq 0}.\label{3.5}
\end{align}
If, instead of \eqref{3.2}, we assume that 
\begin{align}
|\nabla \gamma|\leq A\gamma^{\frac{1}{2}} \ \ \text{outside some compact subset of}\ M, \label{3.6}
\end{align} 
then we can strengthen the second conclusion of \eqref{3.5} to
\begin{align}
\lim_{k\rightarrow\infty} |\nabla{u(x_k)}|=0.
\end{align}
If, instead of \eqref{3.3}, we suppose that 
\begin{gather}
\text{Re} \nabla d\gamma(X, \overline{X})\leq B\gamma^{\frac{1}{2}}G(\gamma^{\frac{1}{2}})^{\frac{1}{2}}|X|^2 \ \text{ outside some compact subset of}\ M\label{3.8}
\end{gather}
for every $X\in T_{1,0}M$, we may strengthen the third conclusion of \eqref{3.5} to 
\begin{gather}
\limsup_{k\rightarrow\infty} {\text{Re} \nabla du(Y_k, \overline{Y_k})\leq 0}
\end{gather}
for every $Y_k\in T_{1,0}M_{x_k}$ with $|Y_k|=1$, where $\text{Re} \nabla du(X, \overline{X})$ is defined as in \eqref{2.19..}.
\end{Thm}
\proof
We define 
\begin{align}
\phi(t)=e^{\int_0^t G(s)^{-\frac{1}{2}}ds}
\end{align}
which is a well-defined, smooth, positive function satisfying $\lim_{t\rightarrow +\infty} \phi(t)=+\infty$. By a simple computation, we have
\begin{align}
\phi'(t)=G(t)^{-\frac{1}{2}}\phi(t),\ \phi''(t)\leq G(t)^{-1}\phi(t),\label{3.11}
\end{align}
and thus,
\begin{align}
\frac{\phi''(t)}{\phi(t)}-\left(\frac{\phi'(t)}{\phi(t)}\right)^2\leq 0.\label{3.12}
\end{align}
Using the assumption \eqref{3.4} $iv)$ and the first formula of \eqref{3.11}, we obtain
\begin{align}
\frac{\phi'(t)}{\phi(t)}\leq C(tG(t^{\frac{1}{2}}))^{-\frac{1}{2}}\label{3.13}
\end{align}
for some positive constant $C$. Now let us consider the function on $M$
\begin{align}
f_k(x)=\frac{u(x)-u(p)+1}{\phi(\gamma(x))^{\frac{1}{k}}}
\end{align}
where $p$ is a fixed point on $M$ and $k\in \mathbb{N}_+$. Since $f_k(p)>0$ and $\limsup_{x\rightarrow\infty}f_k(x)\leq 0$, $f_k$ attains a positive maximum at some point $x_k\in M$ for each $k\in \mathbb{N}_+$. Then by the maximum principle, we obtain
\begin{gather}
(\nabla \log f_k)(x_k)=0,\label{3.15}\\
\text{Re} (\nabla d \log f_k)(X_k, \overline{X_k}) \leq 0\label{3.16}
\end{gather}
for every $X_k\in T_{1,0}M_{x_k}$ with $|X_k|=1$. From \eqref{3.15}, \eqref{3.16}, it follows that 
\begin{align}
\nabla u(x_k)=\frac{1}{k}\left(u(x_k)-u(p)+1\right)\frac{\phi'(\gamma(x_k))}{\phi(\gamma(x_k))}\nabla\gamma(x_k),\label{3.17}
\end{align}
and 
\begin{align}
\text{Re}(\nabla du)(X_k, \overline{X_k})\leq& \frac{1}{k}(u(x_k)-u(p)+1)\Bigg\{\left(\frac{\phi''(\gamma(x_k))}{\phi(\gamma(x_k))}-\left(\frac{\phi'(\gamma(x_k))}{\phi(\gamma(x_k))}\right)^2\right)\notag\Bigg.\\
&\times\bigg.|X_k\cdot\nabla^H\gamma(x_k)|^2+\frac{\phi'(\gamma(x_k))}{\phi(\gamma(x_k))}\text{Re}(\nabla d\gamma)(X_k, \overline{X_k})\notag\bigg.\\
&\Bigg.+\frac{1}{k}\left(\frac{\phi'(\gamma(x_k))}{\phi(\gamma(x_k))}\right)^2 |X_k\cdot\nabla^H\gamma(x_k)|^2\Bigg\}.
\end{align}
Using \eqref{3.12}, 
\begin{align}
\text{Re}(\nabla du)(X_k, \overline{X_k})\leq& \frac{1}{k}(u(x_k)-u(p)+1)\frac{\phi'(\gamma(x_k))}{\phi(\gamma(x_k))}\bigg\{\text{Re}(\nabla d\gamma)(X_k, \overline{X_k})\notag\bigg.\\
&\bigg.+\frac{1}{k}\frac{\phi'(\gamma(x_k))}{\phi(\gamma(x_k))} |X_k\cdot\nabla^H\gamma(x_k)|^2\bigg\}.\label{3.19}
\end{align}
Taking trace in \eqref{3.19}, we get
\begin{align}
\Delta_b u(x_k)\leq \frac{1}{k}(u(x_k)-u(p)+1)\frac{\phi'(\gamma(x_k))}{\phi(\gamma(x_k))}\left\{\Delta_b \gamma(x_k)+\frac{1}{k}\frac{\phi'(\gamma(x_k))}{\phi(\gamma(x_k))} |\nabla^H\gamma(x_k)|^2\right\}.\label{3.20}
\end{align}
We claim that 
\begin{align}
\limsup_{k\rightarrow +\infty} u(x_k)=\sup_M u.
\end{align}
Indeed, if not, there exists $\delta>0$ and some $\hat{x}\in M$ such that 
\begin{align*}
u(\hat{x})>\limsup_{k\rightarrow\infty}u(x_k)+\delta.
\end{align*}
If there is a subsequence $\{x_{k_i}\}$ of $\{x_k\}$, such that $\lim_{i\rightarrow +\infty}\gamma(x_{k_i})=+\infty$, then for $i$ large enough, we have 
\begin{align}
f_k(\hat{x})=\frac{u(\hat{x})-u(p)+1}{\phi(\gamma(\hat{x}))^{\frac{1}{k}}}>\frac{u(x_{k_i})-u(p)+1+\delta}{\phi(\gamma(x_{k_i}))^{\frac{1}{k}}}>f_k(x_{k_i})
\end{align}
contradicting the definition of $x_{k_i}$. If $\{x_k\}$ lie in a compact set, then up to passing to a subsequence, $\lim_{k\rightarrow +\infty}x_k=\bar{x}$, and thus $u(\hat{x})\geq u(\bar{x})+\delta$. On the other hand, since $f_k(x_k)\geq f_k(\hat{x})$ for each $k$, we deduce that 
$$u(\bar{x})-u(p)+1=\lim_{k\rightarrow+\infty}f_k(x_k)\geq \lim_{k\rightarrow+\infty}f_k(\hat{x})=u(\hat{x})-u(p)+1$$
that is,
$$u(\bar{x})\geq u(\hat{x}) $$
which is a contradiction too. Therefore, by choosing a subsequence, we may assume that
\begin{align}
\lim_{k\rightarrow +\infty} u(x_k)=\sup_M u.
\end{align}
If $\{x_k\}$ lie in some compact subset of $M$, then $u$ attains its maximum at some point, hence the results of Theorem \ref{theorem3.1} hold by the classical maximum principle. So we just need to consider the case that there exists a subsequence $\{x_k\}$ with $\lim_{k\rightarrow+\infty}\gamma(x_k)=+\infty$. According to\eqref{3.2}, \eqref{3.13}, \eqref{3.17}, we have, for $k$ large enough, 
\begin{align}
|\nabla^H u(x_k)|\leq\frac{AC}{k}\frac{u(x_k)-u(p)+1}{G(\gamma(x_k)^{\frac{1}{2}})^{\frac{1}{2}}}	 
\end{align}
whose right-hand-side tends to $0$ as $k\rightarrow+\infty$. Hence 
\begin{align}
\lim_{k\rightarrow+\infty}|\nabla^H u(x_k)|=0.
\end{align}
From \eqref{3.2}, \eqref{3.3}, \eqref{3.13} and \eqref{3.20}, we deduce that
\begin{align}
\Delta_b u(x_k)\leq \frac{1}{k}(u(x_k)-u(p)+1)(BC+\frac{AC^2}{kG(\gamma(x_k)^{\frac{1}{2}})^{\frac{1}{2}}})\label{3.26}
\end{align}
which shows that
\begin{align}
\limsup_{k\rightarrow+\infty} \Delta_b u(x_k)\leq 0
\end{align}
proving the results of \eqref{3.5}, since the right-hand-side of \eqref{3.26} tends to $0$ as $k\rightarrow+\infty$. 

If we assume that \eqref{3.6} holds instead of \eqref{3.2}, from \eqref{3.6}, \eqref{3.13}, \eqref{3.17}, we have, for $k$ large enough, 
\begin{align}
|\nabla u(x_k)|\leq\frac{AC}{k}\frac{u(x_k)-u(p)+1}{G(\gamma(x_k)^{\frac{1}{2}})^{\frac{1}{2}}}	 
\end{align}
whose right-hand-side tends to $0$ as $k\rightarrow+\infty$. Hence 
\begin{align}
\lim_{k\rightarrow+\infty}|\nabla u(x_k)|=0.
\end{align}

If, instead of \eqref{3.3}, we suppose \eqref{3.8} holds, by \eqref{3.2}, \eqref{3.8}, \eqref{3.13}, \eqref{3.19} and Cauchy-Schwarz inequality, it is clear that
\begin{align}
\text{Re}(\nabla du)(X_k, \overline{X_k})&\leq \frac{1}{k}(u(x_k)-u(p)+1)(BC+\frac{AC^2}{kG(\gamma(x_k)^{\frac{1}{2}})^{\frac{1}{2}}})
\end{align}
for each $X_k\in T_{1,0}M_{x_k}$ with $|X_k|=1$. Let $k\rightarrow +\infty$, we get 
\begin{align}
\limsup_{k\rightarrow\infty} {\text{Re} (\nabla du)(X_k, \overline{X_k})\leq 0}.
\end{align}
\qed

\begin{Rk}
\begin{enumerate}[(1)]
\item The proof shows that we just need $\gamma$ to be $C^2$ in a neighborhood of $x_k$. Therefore, using the trick of Calabi (cf. \cite{[Ca]}), this theorem remains valid for the important case $\gamma(x)=r^2(x)$, where $r(x)$ is the Riemannian distance with respect to $g_\theta$ from a fixed point $o$ to $x$.\label{remark3.2}
\item Similar generalized maximum principle for Riemannian manifolds was established by \cite{[PRS-2]}.
\end{enumerate}
\end{Rk}

In order to obtain the generalized maximum principle of Yau's type in pseudo-Hermitian geometry, we need a sub-Laplace comparison theorem as follows.

\begin{Lma}\label{lemma3.2}
Let $(M^{2m+1},H,J,\theta)$ be a complete pseudo-Hermitian manifold and $r$ be the Riemannian distance of $g_\theta$ relative to a fixed point $x_0$ in $M$. Assume that
\begin{align}
Ric(\nabla^{1,0}r, \nabla^{0,1}r)\geq -K(r)|\nabla^{1,0}r|^2, \ \ \ \|A\|_{C^1(B(x_0, r))}\leq L(r)
\end{align}
where $\nabla^{1,0}r\in \Gamma^\infty(T_{1,0}M)$ is the unique complex vector field such that $\nabla^Hr=\nabla^{1,0}r+\nabla^{0,1}r$ and $\nabla^{0,1}r=\overline{\nabla^{1,0}r}$, and $K, L$ are two nonnegative $C^1$ functions on $[0, +\infty)$, $B(x_0,r)=\{x\in M:\ r(x)<r \}$. Set
\begin{align}
F(r)=\frac{1}{2m}\left(K(r)+2m^2L^2(r)+\left(2m^2+2\sqrt{2}m^{\frac{3}{2}}-4m\right)L(r)+4\right) \label{3.31}
\end{align}
and assume that it satisfies that $\inf_{[0, +\infty)}\frac{F'}{F^{\frac{3}{2}}}>-\infty$. Then there exists $D>0$ large enough, having set
\begin{align}
h(r)=\frac{1}{D\sqrt{F(0)}}\left(e^{D\int^r_0\sqrt{F(s)}ds}-1\right),\label{3.34}
\end{align}
we have
\begin{align}
\Delta_b r(x)\leq 2m\frac{h'(r(x))}{h(r(x))} \label{3.41}
\end{align}
for any $x$ outside the cut locus of $x_0$. 
\end{Lma}
\proof
Let $e\in T_{x_0}M$ with $|e|=1$, and $\gamma(t)$ be a geodesic with respect to Riemannian connection of $g_\theta$ such that $\gamma'(0)=e$.  if $e\neq\xi$, we choose the unitary frame field $\{\eta_1=\frac{1}{\sqrt{2}|\nabla^Hr|}(\nabla^Hr-\sqrt{-1}J\nabla^Hr), \eta_2,..., \eta_m\}$ of $T_{1,0}M$ on a small neighborhood of $\gamma$. If $e=\xi$, we choose any unitary frame field of $T_{1,0}M$ which is also denoted by $\{\eta_1, \eta_2,..., \eta_m\}$.
Since $|\nabla r|=1$ outside the cut locus of $x_0$, taking the covariant derivatives and summing yield
\begin{align}
0=\frac{1}{2}\sum_{k=1}^m(|\nabla r|^2)_{k\bar{k}}=&\sum_k|r_{0k}|^2+\sum_{i,k}|r_{ik}|^2+\sum_{i,k}|r_{i\bar{k}}|^2\notag\\
&+\sum_{k}r_0r_{0k\bar{k}}+\sum_{i,k}r_{\bar{i}}r_{ik\bar{k}}+\sum_{i,k}r_ir_{\bar{i}k\bar{k}}.\label{3.36.}
\end{align}
Using the commutative formulas for covariant derivatives (cf. \cite{[RYC]} ), it is easy to obtain the following equalities 
\begin{align}
&\sum_kr_0r_{0k\bar{k}}=r_0(\sum_kr_{k\bar{k}})_0+\sum_{k,l}(r_0r_{kl}A^l_{\bar{k}}+r_0r_lA^l_{\bar{k},k}+r_0r_{\bar{l}\bar{k}}A^{\bar{l}}_k+r_0r_{\bar{l}}A^{\bar{l}}_{k,\bar{k}}),\label{3.35.}\\
&\sum_kr_{\bar{i}}r_{ik\bar{k}}=r_{\bar{i}}(\sum_kr_{k\bar{k}})_i+2\sqrt{-1}\sum_kr_{\bar{k}}r_{0k}+\sum_{l}r_lr_{\bar{i}}R_{i\bar{l}}-2\sqrt{-1}\sum_{k,l}r_{\bar{k}}r_{\bar{l}}A^{\bar{l}}_k,\label{3.36}\\
&\sum_kr_ir_{\bar{i}k\bar{k}}=r_i(\sum_kr_{k\bar{k}})_{\bar{i}}-2\sqrt{-1}(m-1)\sum_jr_jr_iA^j_{\bar{i}}-2\sqrt{-1}\sum_kr_{0\bar{k}}r_k, \label{3.37}
\end{align}
where $i,j,k,l=1,2,...,m$. Substituting \eqref{3.35.}-\eqref{3.37} into \eqref{3.36.} yields
\begin{align}
0&=\frac{1}{2}\sum_{k=1}^m((|\nabla r|^2)_{k\bar{k}}+(|\nabla r|^2)_{\bar{k}k})
=2\sum_k|r_{0k}|^2+2\sum_{i,k}|r_{ik}|^2+2\sum_{i,k}|r_{i\bar{k}}|^2\notag\\
&+<\nabla\triangle_b r, \nabla r>-8\text{Im}\{\sum_kr_{\bar{k}}r_{0k}\}+\sum_{i,j}2r_{\bar{i}}r_jR_{i\bar{j}}-4(m-2)\text{Im}\{\sum_{k,l}r_kr_lA^l_{\bar{k}}\}\notag\\&+4\text{Re}\{\sum_{k,l}r_0r_{\bar{l}\bar{k}}A^{\bar{l}}_k\}+2\text{Re}\{\sum_{k,l}r_0r_lA^l_{\bar{k},k}\}+2\text{Re}\{\sum_{k,l}r_0r_{\bar{l}}A^{\bar{l}}_{k,\bar{k}}\}.
\end{align}
Since $Ric(\nabla^{1,0}r, \nabla^{0,1}r)\geq -K(r)|\nabla^{1,0}r|^2$, $\|A\|_{C^1(B(x_0, r))}\leq L(r)$, and using the Cauchy-Schwarz and Young inequalities, it follows that
\begin{align}
0\geq\frac{1}{2m}(\triangle_b r)^2+<\nabla\triangle_b r, \nabla r>-2m F(r)\label{0.23}
\end{align}
where $F(r)$ is defined as \eqref{3.31}.
Set 
\begin{align}
f(t)=\triangle_br(\gamma(t)),
\end{align}
then \eqref{0.23} can be written as
\begin{align}
0\geq \frac{1}{2m}f^2(t)+f'(t)-2m F(t).\label{3.49}
\end{align}
Set
\begin{align}
g(t)=t^ae^{\int^t_0\frac{f(s)}{2m}-\frac{a}{s}ds}
\end{align}
where $a=\frac{2m-|\pi^He|^2}{2m}$.
Since a smooth Riemannian metric is locally Euclidean, we have
\begin{align}
\lim_{t\rightarrow 0} tf(t)=2m-|\pi^H e|^2
\end{align}
which implies the function $g(t)$ is well-defined on $[0,+\infty)$.
A simple computation yieds
\begin{gather}
g(0)=0,\\
g'=\frac{f}{2m}g,\label{3.52}
\end{gather}
and using \eqref{3.49}, we get
\begin{align}
g''\leq F g.\label{3.53}
\end{align}
On the other hand, since $h$ defined in \eqref{3.34} satisfies $h(0)=0$, $h'(0)=1$ and
\begin{align}
h''-Fh\geq\frac{F}{\sqrt{F(0)}}\left[\inf_{[0,+\infty)}\frac{F'}{2F^{\frac{3}{2}}}+D-\frac{1}{D}\right]\geq 0\label{3.49.}
\end{align}
for $D>0$ large enough.
Define
\begin{align}
\Phi(t)=(gh'-g'h)(t).
\end{align}
From \eqref{3.53}, \eqref{3.49.} and $h\geq0$, it follows that 
\begin{align}
\Phi'(t)=gh''-g''h\geq0,
\end{align}
which implies $\Phi(t)\geq \lim_{t\rightarrow0+}\Phi(t)=0$, i.e., $gh'-g'h\geq0$.
Therefore,
\begin{align}
f(t)=2m\frac{g'}{g}(t)\leq 2m\frac{h'}{h}(t).
\end{align}
\qed

In particular, if $K(r)\equiv k_1\geq0$ and $L(r)\equiv k_2\geq0$ in Lemma \ref{lemma3.2} where $k_1, k_2$ are constant, then we have the following corollary which has been proved in \cite{[CDRZ]}.
 \begin{Cor}\label{corollary3.2}
Let $(M^{2m+1},H,J,\theta)$ be a complete pseudo-Hermitian manifold with pseudo-Hermitian Ricci curvature bounded below by $-k_1\leq0$ and $\|A\|_{C^1}$ bounded above by $k_2\geq0$. Let $r$ be the Riemannian distance of $g_\theta$ relative to a fixed point $x_0$. Then for any $x\in M$ which is not on the cut locus of $x_0$, there exists $C=C(m)$ such that
$$\Delta_b r(x)\leq C\left(\frac{1}{r}+\sqrt{1+k_1+k_2+k_2^2}\right).$$
\end{Cor}
 
By Remark 3.2 \eqref{remark3.2}, taking $\gamma(x)=r^2(x)$ and $G(r)=r^2+1$ in Theorem \ref{theorem3.1}, and combining with Corollary \ref{corollary3.2}, we have the generalized maximum principle of Yau's type in pseudo-Hermitian geometry as follows.
\begin{Thm}\label{Thm 3.2}
Let $(M^{2m+1},H,J,\theta)$ be a complete pseudo-Hermitian manifold with pseudo-Hermitian Ricci curvature bounded below and $\|A\|_{C^1}$ bounded above. Then for every $C^2$ function $u:M\rightarrow \mathbb{R}$ which is bounded from above, there exists a sequence $\{x_k\}\subset M$ such that
\begin{align}
\lim_{k\rightarrow\infty}u(x_k)=\sup_M{u},\ \lim_{k\rightarrow\infty}|\nabla{u(x_k)}|=0, \  \limsup_{k\rightarrow\infty} {\Delta_b u(x_k)\leq 0}.
\end{align}
\end{Thm}
 
Furthermore, if we strengthen the curvature condition, the conclusion of Theorem \ref{Thm 3.2} will be stronger. For establishing the generalized maximum principle of Omori's type in pseudo-Hermitian geometry, we need the following Hessian comparison theorem.
\begin{Lma}[cf. \cite{[DRY]}]{\label{lemma3.5}}
Let $(M^{2m+1},H,J,\theta )$ be a complete pseudo-Hermitian manifold with horizontal sectional curvature $K_{\sec }^{H}$ bounded below by $-k _{1}\leq 0$ and $\Vert A\Vert _{C^{1}}$ bounded above by $k_2\geq 0$. Let $r$ be the Riemannian distance of $(M,g_{\theta })$ relative to a fixed point $x_0$. Then there exists a positive constant $D=D(m, k_{1}, k_{2})$ such that
\begin{equation*}
\text{Re} \nabla dr^2(X, \overline{X})\leq D(1+r)|X|^2
\end{equation*}
holds outside the cut locus of $x_0$ for every $X\in T_{1,0}M$.
\end{Lma}

Taking $\gamma(x)=r^2(x)$ and $G(r)=r^2+1$ in Theorem \ref{theorem3.1} again and combining with Lemma \ref{lemma3.5}, we obtain the generalized maximum principle of Omori's type in pseudo-Hermitian geometry as follows.
\begin{Thm}\label{theorem3.6}
Let $(M^{2m+1},H,J,\theta)$ be a complete pseudo-Hermitian manifold with horizontal sectional curvature $K_{\sec }^{H}$ bounded below and $\|A\|_{C^1}$ bounded above. Then for every $C^2$ function $u:M\rightarrow \mathbb{R}$ which is bounded from above, there exists a sequence $\{x_k\}\subset M$ such that
\begin{gather}
\lim_{k\rightarrow\infty}u(x_k)=\sup_M{u},\
\lim_{k\rightarrow\infty} |\nabla{u(x_k)}|=0,\
\limsup_{k\rightarrow\infty} {\text{Re} \nabla du(X_k, \overline{X_k})\leq 0}
\end{gather}
for every $X_k\in T_{1,0}M_{x_k}$ with $|X_k|=1$.
\end{Thm}

\section{Stochastic completeness for pseudo-Hermitian manifolds}\label{section4}
In this section, we will introduce the stochastic completeness for pseudo-Hermitian manifolds, and give the relationship between the generalized maximum principles and stochastic completeness.

Firstly, we construct the heat semigroup associated to $\Delta_b$ on a pseudo-Hermitian manifold $(M^{2m+1},H,J,\theta)$. According to \eqref{2.20..}, $\Delta_b$ is symmetric and non-positive definite on the space $C_0^\infty(M)$. Set 
\begin{align*}
&{}^bW^{1,2}(M)=\{u\in L^2(M,\mu):\ \nabla^H u\in L^2(M,\mu)\},\\
&{}^bW^{1,2}_0(M)=\text{closure\ of}\ C_0^\infty(M) \ \text{in} \ {}^bW^{1,2}(M)\  \text{with\ respect\ to\ } \|\cdot\|_{{}^bW^{1,2}(M)},\\
&{}^bW^{2,2}_0(M)=\{u\in {}^bW^{1,2}_0(M,\mu): \ \Delta_b u\in L^2(M,\mu)\},
\end{align*}
where $\mu$ is the Borel measure given by the volume form $\theta\wedge(d\theta)^m$, $\|\cdot\|_{{}^bW^{1,2}(M)}$ is a norm on ${}^bW^{1,2}(M)$ defined by 
\begin{align}
\|u\|^2_{{}^bW^{1,2}(M)}=\|u\|^2_{L^2(M,\mu)}+\|\nabla^Hu\|^2_{L^2(M,\mu)}.
\end{align}
Similar to Theorem 4.6 in \cite{[Gr]}, the operator $\Delta_b|_{C_0^\infty(M)}$ can be uniquely extended to the operator  $\Delta_b|_{^bW^{2,2}_0(M)}$ which is self-adjoint and densely defined in $L^2(M,\mu)$. This result can be proved by a small modification of the second proof of Theorem 4.6 in \cite{[Gr]}. Then, by the spectral decomposition theorem, $\Delta_b$ can be expressed by $\Delta_b=-\int^\infty_0 \lambda dE_{\lambda}$ in $L^2(M, \mu)$, where $\{E_\lambda\}$ is the spectral resolution. Define the heat semigroup $(P_t)_{t\geq 0}$ by $P_t=-\int^\infty_0 e^{-\lambda t}dE_{\lambda}$, and it is a one-parameter family of bounded operators on  $L^2(M, \mu)$ ($\|P_t\|\leq 1$) which have the following property:
\begin{Thm} \label{theorem4.1}
The solution of the parabolic Cauchy problem
\begin{equation}
\left\{
\begin{aligned}
&\frac{\partial u}{\partial t}-\Delta_b u=0\ \text{in}\ M\times\mathbb{R}^+\\
&u(x,t)\xrightarrow{L^2}f(x)\ \text{as}\ t\rightarrow 0+
\end{aligned}
\right.
\label{4.1..}
\end{equation}
that satisfies $\|u(\cdot, t)\|_{L^2}<\infty$ is unique and given by $u(x,t)=P_t f(x)$. Furthermore, 
from the hypoellipticity of $\frac{\partial}{\partial t}-\Delta_b$, it follows that $u(x,t)=P_t f(x)\in C^\infty(M\times\mathbb{R}^+)$.
\end{Thm}
\begin{Rk}
When $(M,g_\theta)$ is complete, the above theorem has been discussed in \cite{[BW]}, but the completeness is not really necessary. The above theorem  may be deduced from Theorem 4.9 and Theorem 4.10 in \cite{[Gr]}.
\end{Rk}
On a pseudo-Hermitian manifold, there is a positive smooth kernel associated with $P_t$.
\begin{Thm}\label{theorem4.2}
Let $(M^{2m+1},H,J,\theta)$ be a pseudo-Hermitian manifold and $P_t$ the heat semigroup on $L^2$ associated with $\Delta_b$. Then there exists a so-called heat kernel $p_t(x,y)$ satisfying
\begin{enumerate}[(1)]
\item $p_t(x,y)$ is a positive $C^\infty$ function on $\mathbb{R}^+\times M\times M$.
\item $p_t(x,y)=p_t(y,x)$.
\item $\int_M p_t(x,y) d\mu(y)\leq 1$ for all $x\in M$ and $t>0$, such that
\begin{align}
P_t f(x)=\int_M p_t(x,y) f(y) d\mu(y) \label{4.2..}
\end{align}
for any $f\in L^2(M,\mu)$.
\item For any $f\in C^\infty_0(M)$,
\begin{align}
P_tf(x)=\int_M p_t(x,y)f(y)d\mu(y)\xrightarrow{C^\infty(M)} f(x), \ \text{as}\ t\rightarrow 0+,
\end{align}
i.e., $\lim_{t\rightarrow 0+}\|P_tf-f\|_{C^k(K)}=0$ for any compact subset $K$ of $M$ and any $k\in\mathbb{N}$.
\end{enumerate}
\end{Thm}

\begin{Rk}
Although $(1)-(3)$ were proved in \cite{[St]} under the assumption that the Carnot-Carath\'eodory distance is complete, the completeness is not really necessary. The existence and related properties of the heat kernel may be established by an exhaustion method through precompact domains as in the Riemannian case (cf. \cite{[Do]}, \cite{[Gr]}).
\end{Rk}

\proof[Proof of Theorem \ref{theorem4.2}] 
Let us show the conclusion $(4)$. If $f\in C^\infty_0(M)$, then $f\in \text{dom}\ L^k$ for any $k\in \mathbb{N}+$ where $L=\Delta_b|_{{}^bW^{2,2}_0(M)}$, namely,
\begin{align}
\int_{0}^{\infty} \lambda^{2k} d\|E_\lambda f\|^2<+\infty.
\end{align}
According to the spectral theorem, we have
\begin{equation}
L^k f=\int_0^\infty \lambda^k dE_\lambda f
\end{equation}
and 
\begin{equation}
L^kP_t f=\int_0^\infty \lambda^ke^{-t\lambda} dE_\lambda f,
\end{equation}
thus,
\begin{equation}
\|L^k(P_tf-f)\|^2_{L^2(M)}=\int_0^\infty \lambda^k\left(1-e^{-t\lambda}\right)^2 d\|E_\lambda f\|^2\rightarrow 0,\label{4.7.}
\end{equation}
as $t\rightarrow 0+$ for any $k\in \mathbb{N}+$, since the dominated convergence theorem. Set $v=P_tf-f$, by Theorem \ref{theorem4.1}, we have $v\in C^\infty (M\times \mathbb{R}^+)$, and thus $L^k v\in L^2_{loc}(M,\mu)$. Using the interior regularity result for  $\Delta_b$ (cf. Theorem 18 of \cite{[RS]}) yields that for any $\Omega'\subset\subset \Omega''\subset\subset M$, there exists a constant $C>0$ such that
\begin{align}
\|L^{k-i}v\|_{S^2_{2i}(\Omega')}\leq C\left(\|L^{k-i}v\|_{L^2(\Omega'')}+\|L^{k-i+1}v\|_{S^2_{2(i-1)}(\Omega'')} \right)
\end{align}
where $i=1,2,\dots, k$. Therefore, 
\begin{align}
\|v\|_{S^2_{2k}(\Omega')}\leq C\sum_{l=0}^k \|L^lv\|_{L^2(\Omega'')}
\end{align}
for any positive integer $k$. By the embedding theorem $S^2_{2k}(\Omega')\subset W^{k,2}(\Omega')\subset C^{s}(\Omega')$ ($k>s+\frac{2m+1}{2}$) (cf. Theorem 19.1 of \cite{[FS]}) where $W^{k,2}(\Omega')$ is the classical Sobolev space, we have
\begin{align}
\|v\|_{C^{s}(\Omega')}\leq C\sum_{l=0}^k \|L^lv\|_{L^2(\Omega'')}\rightarrow 0
\end{align}
as $t\rightarrow 0+$, due to Theorem \ref{theorem4.1} and  \eqref{4.7.}. Hence, $v=P_tf-f\xrightarrow{C^\infty}0$ as $t\rightarrow 0+$.
\qed

Using the identity \eqref{4.2..}, we can extend the definition of the operators $P_t$ as follows. Set 
\begin{align}
P_t f(x)=\int_M p_t(x,y) f(y) d\mu(y) 
\end{align}
for any function $f$ such that the right hand side of the above identity makes sense. So the operators $P_t$ transform positive functions into positive functions and satisfy $0<P_t 1\leq 1$. Furthermore, we can solve the Cauchy problem \eqref{4.1..} when $f$ is a bounded continuous function on $M$ by using the heat semigroup $P_t$ whose proof is similar to Theorem 7.16 and Corollary 8.3 of \cite{[Gr]}.
\begin{Thm}\label{theorem4.4.}
For any bounded continuous function $f$ on $M$, $u(x,t)=P_tf(x)$ is a bounded $C^\infty(M\times\mathbb{R}^+)$ solution of 
\begin{equation}
\left\{
\begin{aligned}
&\frac{\partial u}{\partial t}-\Delta_b u=0, \ \ \text{in}\ M\times\mathbb{R}^+,\\
&u(x,t)\rightarrow f(x),\ \ \text{as}\ t\rightarrow 0+,
\end{aligned}
\right.
\label{4.13.}
\end{equation}
where the initial condition means that $\lim_{t\rightarrow 0+}\|u(x,t)-f(x)\|_{C^0(K)}=0$ for any compact set $K\subset M$.
If, in addition, $f\geq 0$, then $u(t,x)=P_tf(x)$ is its minimal nonnegative solution.
\end{Thm}

Analogue to the stochastic completeness in Riemannian Geometry, one can also introduce a similar definition for pseudo-Hermitian manifolds.
\begin{Def}
 A pseudo-Hermitian manifold $(M,H,J,\theta)$ is said to be stochastically complete if $P_t 1=1$ for all $t>0$.
\end{Def}
\begin{Rk}\label{remark4.6}
By the semigroup property of $P_t$, one can show that $P_t 1=1$ for all $t>0$ if and only if $P_t 1(x)=1$ for some $(x,t)\in M\times\mathbb{R}^+$. (cf. Theorem 6.2 of \cite{[Gr-1]})
\end{Rk}

Now we discuss the relationship between stochastic completeness and generalized maximal principles. For this purpose, we need the following lemmas.

\begin{Lma}\label{lemma4.5}
Let $\Omega$ be a relatively compact connected open subset of a pseudo-Hermitian manifold $(M^{2m+1},H,J,\theta)$. Let $Lu=\Delta_b u-\lambda u$, where $\lambda$ is a positive constant. If $u$ is a ${}^bW^{1,2}(\Omega)$ function satisfying 
\begin{equation}
Lu\geq 0,\ \text{in}\ \Omega,\label{4.14..}
\end{equation}
\begin{equation}
u\leq 0\mod{}^bW^{1,2}_0(\Omega), \label{4.2.}
\end{equation}
where \eqref{4.14..} is understood in the sense of distributions and the boundary condition \eqref{4.2.} means that $u\leq u_0$ for some $u_0\in{}^bW^{1,2}_0(\Omega)$, then $u\leq 0$ in $\Omega$.
\end{Lma}
\proof By a small modification of the proof of Lemma 5.12 in \cite{[Gr]}, it is easy to show that \eqref{4.2.} holds if and only if $u_+=\max{(u,0)}\in {}^bW^{1,2}_0(\Omega)$. Hence, taking $u_+$ as a test function of \eqref{4.14..} yields 
\begin{align}
-\int_\Omega L_\theta(\nabla^Hu, \nabla^Hu_+)d\mu-\lambda\int_\Omega uu_+d\mu\geq0.
\end{align}
Since 
\begin{align}
\nabla^Hu_+=
\begin{cases}
\nabla^H u & \mbox{if }u>0,\\
0 & \mbox{if }u\leq0,
\end{cases}
\end{align}
we have
\begin{align}
\int_\Omega |\nabla^Hu_+|^2d\mu\leq -\lambda\int_\Omega (u_+)^2d\mu\leq 0.\label{4.18..}
\end{align}
On the other hand, by the Poincar\'e inequality (cf. \cite{[JX]}), there is a constant $C>0$ such that
\begin{align}
\int_\Omega |\nabla^Hu_+|^2d\mu\geq C\int_\Omega (u_+)^2d\mu.\label{4.19..}
\end{align}
Combining \eqref{4.18..} and \eqref{4.19..} yields that $u_+=0$, i.e., $u\leq 0$.
\qed
\begin{Lma}\label{lemma4.6}
Let $(M^{2m+1}, H, J, \theta)$ be a non-compact pseudo-Hermitian manifold. Assume that $H(x,u)$ is a continuous function on $M\times \mathbb{R}$ and is locally Lipschitz continuous with respect to $u$ when $x$ remains in a compact subset of $M$. If $u_+$ and $u_-$ are two $C^0(M)\cap {}^bW^{1,2}_{loc}(M)$ functions satisfying 
\begin{equation}
\Delta_b u_++H(x, u_+)\leq 0\ \ \ \ \text{in}\ M,\label{4.9..}
\end{equation}
\begin{equation}
\Delta_b u_-+H(x, u_-)\geq 0\ \ \ \ \text{in}\ M,\label{4.6}\\
\end{equation}
\begin{equation}
u_+\geq u_- \ \ \ \ \text{in}\ M,\label{4.7}
\end{equation}
then there exists a function $u\in C^0(M)$ such that
\begin{equation}
\Delta_b u+H(x, u)=0\ \ \ \ \text{in}\ M,\label{4.8.}
\end{equation}
\begin{equation}
u_-\leq u\leq u_+ \ \ \ \  \text{in}\ M.\label{4.9.}
\end{equation}
Here \eqref{4.9..}, \eqref{4.6}, \eqref{4.8.} are understood in the sense of distributions.
\end{Lma}
\proof
Since $M$ is connected and has a countable topological basis, there exists a sequence $\{\Omega_k\}_{k=1}^\infty$ of relatively compact connected open subset of $M$ such that $\Omega_k\subset\subset \Omega_{k+1}$ and $M=\bigcup_{k}\Omega_k$.
We firstly construct a solution in $\Omega_n$ by monotone iteration schemes (cf. \cite{[Ni]}, \cite{[RRV]}, \cite{[Sa]}). Set $r_n^-=\min_{\Omega_n} u_-$, $r_n^+=\max_{\Omega_n} u_+$, $I_n=[r_n^-, r_n^+]$. Since $H(x,r)$ is locally Lipschitz continuous with respect to $r$ when $x\in \Omega_n$, there exists $\lambda_n\in\mathbb{R}^+$ and $H_n(x, r)\in C^0(M\times \mathbb{R})$ such that $H(x, r)=H_n(x, r)-\lambda_n r$ and the function $H_n(x,r)$ is increasing with respect to $r\in I_n$ when $x\in \Omega_n$. We want to look for a function $u_n$ satisfying
\begin{equation}
\Delta_b u_n -\lambda_n u_n=-\tilde{H}_n(x, u_n) \ \ \ \text{in}\ \Omega_n\label{4.11}
\end{equation}
\begin{equation}
u_-(x)\leq u_n(x)\leq u_+(x)\ \ \ \text{in}\ \Omega_n\label{4.12.}
\end{equation}
where 
\begin{equation}
\tilde{H}_n(x, r)=
\begin{cases}
H_n(x, r_n^-) &\mbox{if} \ r<r_n^-,\\
H_n(x, r) & \mbox{if}\  r\in I_n,\\
H_n(x, r_n^+)  & \mbox{if}\  r>r_n^+.
\end{cases}
\end{equation}
In order to find this function, we consider the sequences $\{v_-^k\}\subset {}^bW^{1,2}(\Omega_n)$ defined for $k\geq1$ by  
\begin{equation}
\left\{ 
\begin{aligned}
&(\Delta_b-\lambda_n id) v_-^k=-\tilde{H}_n(x, v_-^{k-1}) \ \ \ \text{in}\ \Omega_n\\
& v_-^k-\frac{u_-+u_+}{2}\in{}^bW^{1,2}_0(\Omega_n)\\
&v_-^0=u_-.
\end{aligned}
\right.
\label{4.18}
\end{equation}
From \eqref{4.9..} and \eqref{4.18}, we deduce that
\begin{equation}
\left\{ 
\begin{aligned}
&(\Delta_b-\lambda_n id) (v_-^1-u_+)\geq \tilde{H}_n(x, u_+)-\tilde{H}_n(x, u_-)\geq0 \ \ \ \text{in}\ \Omega_n\\
&v_-^1-u_+ \leq 0\ \mod{}^bW^{1,2}_0(\Omega_n). \\
\end{aligned}
\right.
\label{}
\end{equation}
According to Lemma \ref{lemma4.5}, we have 
\begin{align}
v_-^1\leq u_+.
\end{align}
By \eqref{4.6} and \eqref{4.18}, we get 
\begin{equation}
\left\{ 
\begin{aligned}
&(\Delta_b-\lambda_n id) (u_--v_-^1)=\left(\Delta_b u_-+H(x, u_-)\right)\geq 0 \ \ \ \text{in}\ \Omega_n\\
&u_--v_-^1 \leq 0\ \mod{}^bW^{1,2}_0(\Omega_n). \\
\end{aligned}
\right.
\label{}
\end{equation}
Using Lemma \ref{lemma4.5} again, we have
\begin{align}
u_-\leq v_-^1.
\end{align}
Consequently, 
\begin{align}
u_-\leq v_-^1\leq u_+.
\end{align}
Iterating the above proceeding gives
\begin{align}
u_-\leq v_-^1\leq v_-^2\leq \ldots\leq u_+.
\end{align}
Set $\underline{u}_n=\lim_{k\rightarrow +\infty} v_-^k$. Then by the dominated convergence theorem, we get that $\underline{u}_n$ is a solution of \eqref{4.11} and \eqref{4.12.} in the distributional sense. Hence, $\underline{u}_n$ solves in the distributional sense the problem
\begin{equation}
\left\{
\begin{aligned}
&\Delta_b u-H(x, u)=0\ \ \ \ \text{in}\ \Omega_n,\\
&u_-\leq u\leq u_+ \ \ \ \  \text{in}\ \Omega_n.
\end{aligned}
\right.
\label{4.20.}
\end{equation}
 Since $\underline{u}_n\in I_2$ on $\Omega_2$ for $n\geq 2$ and $H_n(x, r)\in C^0(M\times \mathbb{R})$, there exists a constant $C>0$ such that $\|H(x,\underline{u}_n(x))\|_{L^\infty(\Omega_2)}\leq C$ and $\| \underline{u}_n\|_{L^\infty(\Omega_2)}\leq C$ for all $n\geq 2$. Using the interior regularity result for the sub-Laplacian operator (cf. Theorem 18 of \cite{[RS]}), we get that for any $p>1$ and any relatively compact open set $\Omega'$ with $\Omega_1\subset \subset\Omega'\subset \subset\Omega_2$, there exists $C_p>0$ such that
 \begin{align}
 \|\underline{u}_n\|_{S_2^p(\Omega')}\leq C_p.
 \end{align}
Taking $p>2m+1$ and using the embedding theorem $S^p_2(\Omega')\subset W^{1,p}(\Omega')\subset\subset C^0(\Omega')$ (cf. Theorem 19.1 of \cite{[FS]}), where $W^{1,p}(\Omega')$ is the classical Sobolev space, we obtain that there is a function $u_1\in{C^0(\Omega_1)}$ such that
\begin{align}
\|\underline{u}_{n_{1j}}-u_1\|_{C^0(\Omega_1)}\rightarrow 0
\end{align}
as $j\rightarrow \infty$, where $\underline{u}_{n_{1j}}$ is a subsequence of $\underline{u}_n$. Then $u_1$ is a $C^0$-solution of \eqref{4.20.} for $n=1$. By Repeating the above argument, there is a subsequence  $\underline{u}_{n_{2j}}$ of  $\underline{u}_{n_{1j}}$ converging to $u_2$ in $C^0(\Omega_2)$, and $u_2$ is a $C^0$-solution of \eqref{4.20.} for $n=2$ with $u_2|_{\Omega_1}=u_1$ .  In this way, we get a subsequence $\underline{u}_{n_{kj}}$ converging to $u_k$ in $C^0(\Omega_k)$, and $u_k$ is a $C^0$-solution of \eqref{4.20.} for $n=k$ with $u_k|_{\Omega_{k-1}}=u_{k-1}$. By a diagonal process, $\{\underline{u}_{n_{mm}}\}_{m=k}^\infty$ is a subsequence of $\{\underline{u}_{n_{kj}}\}_{j=1}^\infty$ for every $k$. Thus, we have $\lim_{m\rightarrow\infty}\underline{u}_{n_{mm}}=u_k$ on $\Omega_k$ for any $k$. 
Therefore, $u(x)=\lim_{m\rightarrow\infty}\underline{u}_{n_{mm}}(x)$ is a $C^0$-solution of \eqref{4.8.} and \eqref{4.9.}.
\qed

In particular, taking $H(x,u)=\lambda u\  (\lambda\in \mathbb{R})$ in the above lemma yields

\begin{Cor}\label{corollary4.7}
Let $(M^{2m+1}, H, J, \theta)$ be a non-compact pseudo-Hermitian manifold.
If $u_+$ and $u_-$ are two $C^0(M)\cap {}^bW^{1,2}_{loc}(M)$ functions satisfying
\begin{equation}
\Delta_b u_+-\lambda u_+\leq 0\ \ \ \ \text{in}\ M,\label{4.38..}
\end{equation}
\begin{equation}
\Delta_b u_--\lambda u_-\geq 0\ \ \ \ \text{in}\ M,\label{4.39..}\\
\end{equation}
\begin{equation}
u_+\geq u_- \ \ \ \ \text{in}\ M,
\end{equation}
where \eqref{4.38..}, \eqref{4.39..} are understood in the sense of distribution, constant $\lambda\in \mathbb{R}$,   then there exists a function $u\in C^\infty(M)$ such that
\begin{equation}
\Delta_b u-\lambda u=0\ \ \ \ \text{in}\ M,
\end{equation}
\begin{equation}
u_-\leq u\leq u_+ \ \ \ \  \text{in}\ M.
\end{equation}
\end{Cor}
\proof
By Lemma \ref{lemma4.6}, we obtain that there exists a function $u\in C^0(M)$ such that $\Delta_b u-\lambda u=0$ and $ u_-\leq u\leq u_+$ in $M$, so it suffices to show $u\in C^\infty(M)$. Let $x$ be an arbitrary point in $M$ and $\Omega$ a relatively compact open neighborhood of $x$ in $M$.
Since $u\in C^0(M)$, $\Delta_b u=\lambda u\in L^p(\Omega)$. From the interior regularity result for the sub-Laplacian operator (cf. Theorem 18 of \cite{[RS]}), it follows that $u\in S^p_2(\Omega_1)$ where $x\in\Omega_1\subset\subset \Omega$. Repeating the above argument, we get that $u\in S^p_{2k+2}(\Omega_{k+1})$ where $x\in\Omega_{k+1}\subset\subset \Omega_{k}\subset\subset\ldots \subset\subset\Omega_1\subset\subset \Omega$ and any $k\geq 1$. According to $S^p_{2k+2}(\Omega_{k+1})\subset C^{k}(\Omega_{k+1})\ (p>2m+1)$ (cf. Theorem 19.1 of \cite{[FS]}), $u\in C^k(\Omega_{k+1})$ for $k\geq 1$, and thus $u$ is $C^\infty$ at $x$. For the arbitrariness of $x$, $u$ is smooth in $M$.
\qed

Now we can give the main theorem of this section, which is a pseudo-Hermitian counterpart of Theorem 1.1 in \cite{[PRS-1]}.

\begin{Thm}\label{theorem4.9}
Let $(M^{2m+1},H,J,\theta)$ be a non-compact pseudo-Hermitian manifold. The following statements are equivalent:
\begin{enumerate}[(1)]
\item $M$ is stochastically complete.
\item For every $\lambda>0$, the only nonnegative bounded $C^2$ function on $M$ with $\Delta_b  u=\lambda u$ is $u\equiv0$.
\item If $u$ is a $C^2$ function on $M$ with $\sup_M u<+\infty$, then $\inf_{\Omega_\alpha} \Delta_b u\leq 0$ for any $\alpha>0$, where $\Omega_\alpha=\{x\in M:\ u(x)>\sup_M u-\alpha \}$.
\item For any $u\in C^2(M)$ with $\sup_M u<+\infty$, there exists $\{x_n\}\subset M$ such that $u(x_n)\geq \sup_M u-\frac{1}{n}$ and $\Delta_b u(x_n)\leq \frac{1}{n}$.
\end{enumerate}
\end{Thm}
\proof $(1)\Rightarrow (2):$ We prove by contradiction, and assume that $(2)$ does not hold, that is, there exists a number $\lambda>0$ such that the equation $\Delta_b  \nu=\lambda \nu$ admits a nontrivial nonnegative bounded $C^2$ solution. Then $u(x,t)=\nu(x)e^{\lambda t}$ solves the following problem:
\begin{equation}
\left\{ 
\begin{aligned}
&\frac{\partial u}{\partial t}=\Delta_b u\ \ \ \text{ on}\  M\times (0, T)\\
&u(x,t)\rightarrow \nu(x)\ \text{as}\ t\rightarrow 0+.
\end{aligned}
\right.
\end{equation}
On the other hand, according to Theorem \ref{theorem4.4.}, there is another solution $\omega=P_t\nu$ to the above problem. Since for any $t>0$,
\begin{align}
\sup_{x\in M} |\omega(x,t)|\leq \sup_{x\in M}|\nu(x)|\int_M p_t(x,y)d\mu(x)=\sup_{x\in M}|\nu(x)| P_t1=\sup_{x\in M}|\nu(x)|
\end{align}
and 
\begin{align}
\sup_{x\in M} |u(x,t)|=e^{\lambda t} \sup_{x\in M}|\nu(x)|>\sup_{x\in M}|\nu(x)|,
\end{align}
then $u(\cdot, t)\not\equiv\omega(\cdot, t)$ for any $t>0$, so the function $f=u-\omega$ is a nontrivial solution of  
\begin{equation}
\left\{ 
\begin{aligned}
&\frac{\partial f}{\partial t}=\Delta_b f \ \ \ \text{ on}\  M\times (0, T)\\
&f(x,t)\rightarrow0\ \text{as}\ t\rightarrow 0+
\end{aligned}
\right.\label{4.4.}
\end{equation}
where $0<T<+\infty$. Since $f$ is non-zero and bounded on $M\times (0, T)$, we assume that $\sup_{M\times (0, T)}f>0$ and $\sup_{M\times (0, T)}|f|<1$. Set $g=1-f$, so $g>0$ and $\inf_{M\times (0,T)} g<1$. According to \eqref{4.4.}, it is easy to see that $g\not\equiv 1$ solves 
\begin{equation}
\left\{ 
\begin{aligned}
&\frac{\partial g}{\partial t}=\Delta_b g\ \ \ \text{ on}\  M\times (0, T)\\
&g(x,t)\rightarrow1\ \text{as}\ t\rightarrow 0+.
\end{aligned}
\right.\label{4.5.}
\end{equation}
By Theorem \ref{theorem4.4.}, we obtain that $P_t 1$ is the minimal positive solution of \eqref{4.5.}, so $P_t 1\leq g$. From  $\inf_{M\times (0,T)} g<1$, we have there exists some $x_0\in M$ and $t_0\in (0, T)$ such that $g(x_0,t_0)<1$, hence $P_{t_0}(x_0)<1$ which contradicts the stochastic completeness of $M$.

$(2)\Rightarrow (1):$ Assume that $M$ is not stochastically complete, namely, $0<P_t1<1$ for all $t>0$ according to Remark \ref{remark4.6}. Set $u(x,t)=P_t 1$ and $v=1-\lambda\int^{+\infty}_0 e^{-\lambda t }u(x,t)dt$ where $\lambda>0$. Then
\begin{align}
\Delta_b v&=-\lambda\int^{+\infty}_0 e^{-\lambda t }\Delta_b u(x,t)dt\notag\\
                &=-\lambda\int^{+\infty}_0 e^{-\lambda t }\frac{\partial }{\partial t} u(x,t)dt\notag\\
                &=-\lambda e^{-\lambda t } u(x,t)|^{+\infty}_0-\lambda^2 \int^{+\infty}_0 e^{-\lambda t }u(x,t)dt\notag\\
                &=\lambda v.
\end{align}
By $0<u(x,t)=P_t 1<1$, we have 
\begin{align}
0<\int^{+\infty}_0 e^{-\lambda t }u(x,t)dt<\int^{+\infty}_0 e^{-\lambda t }dt=\lambda^{-1},
\end{align}
Thus, $0<v<1$ which leads to a contradiction with the assumption $(2)$ of Theorem \ref{theorem4.9}.

$(2)\Rightarrow (3):$ Assume that $(3)$ doesn't hold, that is, there exists a function $u\in C^2(M)$ such that $\sup_M u<+\infty$ and $\inf_{\Omega_\alpha} \Delta_b u\geq 2c>0$ for some $\alpha>0$. Let $\hat{\Omega}=\{x\in M:\ \Delta_b u>c \}$, so $\bar{\Omega}_\alpha\subset \hat{\Omega}$. Set $v=u-\sup_M u+\alpha$. It is obvious that
\begin{align}
\Delta_b v=\Delta_b u>\frac{c}{\alpha}v
\end{align}
on $\hat{\Omega}$, since $v\leq \alpha$. Hence, $v$ is a $C^2$ subsolution of $\Delta_b u=\frac{c}{\alpha}u$ on $\hat{\Omega}$. Since the function $0$ is a subsolution of $\Delta_b u=\frac{c}{\alpha}u$ on $M$, we see that $u_\alpha=\max\{v, 0\}=\max\{u-\sup_M u+\alpha, 0\}\in C^0(M)\cap {}^bW^{1,2}_{loc}(M)$ is also a subsolution on $M$. Note that $u_\alpha\not\equiv0$ and $0\leq u_\alpha\leq \alpha<+\infty$. Then choosing a positive constant $\beta>\alpha$, which is a supersolution of $\Delta_b u=\frac{c}{\alpha}u$ on $M$, and applying Corollary \ref{corollary4.7} yields a smooth solution $w$ on $M$ with $0\leq u_\alpha\leq w\leq \beta$. Since $u_\alpha\not\equiv 0$, so does $w$, which is a contradiction with $(2)$.

$(3)\Rightarrow (4):$ Let $\alpha=\frac{1}{n}$. From $\inf_{\Omega_{\frac{1}{n}}} \Delta_b u\leq 0$ for any $n\in\mathbb{N}$, it follows that any for $n\in \mathbb{N}$, there exists $x_n\in \Omega_{\frac{1}{n}}$, namely, $x_n\in M$ satisfies $u(x_n)\geq \sup_M u-\frac{1}{n}$, such that $\Delta_b u(x_n)\leq \frac{1}{n}$.

$(4)\Rightarrow (2):$ We assume that $u\in C^2(M)$ satisfies $\Delta_b  u=\lambda u$ and $0\leq u\leq \sup_M u<+\infty$. By (4), we get that there exists $\{x_n\}\subset M$ such that $u(x_n)\geq \sup_M u-\frac{1}{n}$ and $\Delta_b u(x_n)\leq \frac{1}{n}$.  Therefore, $\lambda\sup_M u-\frac{\lambda}{n}\leq \lambda u(x_n)=\Delta_b u(x_n)\leq \frac{1}{n}$, and then let $n\rightarrow +\infty$, we have $ \sup_M u\leq 0$, and thus, $u\equiv 0$.
\qed

Combining with the Theorem \ref{theorem3.1}, we have 
\begin{Cor}
Let $(M^{2m+1},H,J,\theta)$ be a pseudo-Hermitian manifold with the hypotheses \eqref{3.1}, \eqref{3.2}, \eqref{3.3}, \eqref{3.4}. Then $M$ is stochastically complete.
\end{Cor}
According to Theorem \ref{Thm 3.2} and Theorem \ref{theorem3.6}, we obtain 
\begin{Cor}
Let $(M^{2m+1},H,J,\theta)$ be a complete pseudo-Hermitian manifold with pseudo-Hermitian Ricci curvature bounded below or horizontal sectional curvature $K_{\sec }^{H}$ bounded below, and $\|A\|_{C^1}$ bounded above. Then $M$ is stochastically complete.
\end{Cor}

Before ending this section, we would like to mention that the stochastically completeness of the heat semigroup $P_t$ generated by the sub-Laplacian on a contact manifold was investigated in \cite{[BW]} too. They gave a sufficient condition, in terms of a generalized curvature inequality, to deduce the stochastically completeness.  Finally it is also obvious that our Theorem \ref{theorem4.9} may be generalized to more general sub-Riemannian manifolds.

\section{Application}
In this section, we consider some applications of the generalized maximal principles in pseudo-Hermitian geometry. Firstly, we use them to study differential inequalities on a pseudo-Hermitian manifold.
\begin{Thm}\label{theorem5.1}
Let $(M^{2m+1},H,J,\theta)$ be a pseudo-Hermitian manifold with hypotheses \eqref{3.1}, \eqref{3.3}, \eqref{3.4}, \eqref{3.6}. Assume that $u$ is a $C^2$ function on $M$ satisfying
\begin{align}
\triangle_b u\geq \phi(u, |\nabla^L u|, |\nabla^H u|),\label{5.1}
\end{align}
where $\nabla^Lu=\pi_L(\nabla u)$, $\phi(t,x_1,x_2)$ is $C^0$ in $t$, $C^2$ in $(x_1,x_2)$, and such that $(\text{Hess}_{x_1,x_2}\phi)(t,x_1,x_2)\geq 0$. Set $f(t)=\phi(t,0,0)$. Let $F$ be a continuous function which is positive on some interval $[a, +\infty)$ and satisfies 
\begin{equation}
\int_{b}^{\infty}\left\{\int_a^t F(s)ds\right\}^{-\frac{1}{2}}dt<+\infty \ \text{for\ some}\ b>a,\label{5.2}
\end{equation}
\begin{equation}
\limsup_{t\rightarrow +\infty} \frac{\int_a^t F(s)ds}{tF(t)}<+\infty,\label{5.3}
\end{equation}
\begin{equation}
\liminf_{t\rightarrow +\infty} \frac{\{\int_a^t F(s)ds\}^{\frac{1}{2}}}{F(t)}\frac{\partial \phi}{\partial x_i}(t,0,0)>-\infty,\  i=1,2.\label{5.4}
\end{equation}
\begin{enumerate}[(1)]
\item If 
\begin{equation}
\liminf_{t\rightarrow +\infty} \frac{f(t)}{F(t)}>0,\label{5.5}
\end{equation}
then $u$ is bounded above on $M$ and $f(\sup_M u)\leq 0$.
\item If
\begin{equation}
\liminf_{t\rightarrow +\infty} \frac{f(t)}{F(t)}\leq0,
\end{equation}
then  either $\sup_Mu=+\infty$ or $u$ is bounded above on $M$ and $f(\sup_M u)\leq 0$.
\end{enumerate}
\end{Thm}
\proof
Consider the following function on $M$:
\begin{align}
\psi=\frac{1}{g(u)}
\end{align}
where 
\begin{equation}
g(t)=
\begin{cases}
\int_{b}^{t}\left\{\int_a^{s} F(r)dr\right\}^{-\frac{1}{2}}ds+2 &\mbox{on}\ [b,+\infty)\\
\mbox{increasing from 1 to 2} &\mbox{on}\ (-\infty, b]
\end{cases}
\end{equation}
with $g\in C^2(\mathbb{R})$ and $g''<0$ on $\mathbb{R}$. Note that $g\geq 1$ and $g'\geq 0$.  Since $\psi$ is bounded blew on M, by Theorem \ref{theorem3.1}, there exists a sequence $\{x_n\}\subset M$ such that
\begin{align}
&\lim_{n\rightarrow +\infty} \psi(x_n)=\inf_M \psi \label{5.9},\\
&|\nabla^L \psi|(x_n)=\frac{g'(u)|\nabla^L u|}{g^2(u)}(x_n)\leq \frac{1}{n}\label{5.10},\\
&|\nabla^H \psi|(x_n)=\frac{g'(u)|\nabla^H u|}{g^2(u)}(x_n)\leq \frac{1}{n}\label{5.11},\\
&\Delta_b \psi(x_n)=-\frac{g''(u)|\nabla^H u|^2}{g^2(u)}-\frac{g'(u)\Delta_b u}{g^2(u)}+\frac{2\left(g'(u) \right)^2|\nabla^H u|^2}{g^3(u)}\geq -\frac{1}{n}.\label{5.12}
\end{align}
Multiplying \eqref{5.12} by $\frac{\left[g'(u)\right]^2}{g^2(u)|g''(u)|}(x_n)$ and substituting \eqref{5.1}, \eqref{5.11} into it, we obtain at $x_n$
\begin{align}
\frac{\left[g'(u)\right]^3\phi(u,|\nabla^L u|, |\nabla^H u|)}{g^4(u)|g''(u)|}\leq \frac{\left[g'(u)\right]^2}{g(u)|g''(u)|}\frac{1}{n}\left(1+\frac{2}{n}\right)+\frac{1}{n^2}.\label{5.13}
\end{align}
By Taylor's formula for $\phi(t,x_1,x_2)$ with respect to $(x_1,x_2)$ and $(Hess_{x_1,x_2}\phi)\geq 0$,
\begin{align}
\phi(t,|\nabla^L u|,|\nabla^H u|)\geq f(u)+\frac{\partial \phi}{\partial x_1}(u,0,0)|\nabla^L u|+\frac{\partial \phi}{\partial x_2}(u,0,0)|\nabla^H u|.
\end{align}
Thus the left hand side of \eqref{5.13} is bounded blew by
\begin{align}
\frac{\left[g'(u)\right]^3f(u)}{g^4(u)|g''(u)|}+\frac{\left[g'(u)\right]^3}{g^4(u)|g''(u)|}\frac{\partial \phi}{\partial x_1}(u,0,0)|\nabla^L u|+\frac{\left[g'(u)\right]^3}{g^4(u)|g''(u)|}\frac{\partial \phi}{\partial x_2}(u,0,0)|\nabla^H u|.
\end{align}
By \eqref{5.10} and \eqref{5.11}, we have
\begin{align}
\frac{\left[g'(u)\right]^3}{g^4(u)|g''(u)|}\frac{\partial \phi}{\partial x_1}(u,0,0)|\nabla^L u|\geq \min\left\{0, \frac{1}{n}\frac{\left[g'(u)\right]^2}{g^2(u)|g''(u)|}\frac{\partial \phi}{\partial x_1}(u,0,0) \right\}=:\frac{1}{n}w_1(u),
\end{align}
and
\begin{align}
\frac{\left[g'(u)\right]^3}{g^4(u)|g''(u)|}\frac{\partial \phi}{\partial x_1}(u,0,0)|\nabla^H u|\geq \min\left\{0, \frac{1}{n}\frac{\left[g'(u)\right]^2}{g^2(u)|g''(u)|}\frac{\partial \phi}{\partial x_2}(u,0,0) \right\}=:\frac{1}{n}w_2(u).
\end{align}
Therefore, at $x_n$, 
\begin{align}
\frac{\left[g'(u)\right]^3f(u)}{g^4(u)|g''(u)|}+\frac{1}{n}w_1(u)+\frac{1}{n}w_2(u)\leq \frac{\left[g'(u)\right]^2}{g(u)|g''(u)|}\frac{1}{n}\left(1+\frac{2}{n}\right)+\frac{1}{n^2}.\label{5.18}
\end{align}

To prove conclusion (1), we argue by contradiction. Assume that $\sup_M u=+\infty$. Since $g'\geq 0$, \eqref{5.9} gives that $\lim_{n} u(x_n)=\sup_M u$, so $\lim_n r(x_n)=+\infty$ due to the continuity of $u$, where $r$ is the Riemannian distance with respect to the Webster metric. 
According to \eqref{5.3} and $g(t)\geq (t-b)\left\{\int^t_a F(r)dr\right\}^{-\frac{1}{2}}$ for $t>b$,
\begin{align}
\limsup_{t\rightarrow +\infty} \frac{\left[g'(t)\right]^2}{g(t)|g''(t)|}\leq \limsup_{t\rightarrow +\infty}\frac{2\left\{\int^t_a F(r)dr\right\}^{\frac{1}{2}}}{(t-b)F(t)}<+\infty.
\end{align}
By \eqref{5.2} and \eqref{5.5},
\begin{align}
\liminf_{t\rightarrow +\infty}\frac{\left[g'(t)\right]^3f(t)}{g^4(t)|g''(t)|}=\liminf_{t\rightarrow +\infty} \frac{2f(t)}{g^4(t)F(t)}>0.
\end{align}
Moreover, using \eqref{5.2} and \eqref{5.4} yields
\begin{align}
w_i(u(x_n))\geq -B^2_i, 
\end{align}
where $B_i$ is independent of $n$ for $i=1,2$. Let $n\rightarrow +\infty$ in \eqref{5.18}, we get a desired contradiction. Thus, $\sup_M<+\infty$. Letting $n\rightarrow +\infty$ in \eqref{5.18} again gives $f(\sup_M u)\leq 0$ which, indeed, completes the proof of (1) and (2).

\qed
\begin{Rk}
\begin{enumerate}[(1)]
\item If $\frac{\partial \phi}{\partial x_1}(t,0,0)=0$, the hypothesis \eqref{3.6} in Theorem \ref{theorem5.1} can be weaken to \eqref{3.2}.
\item Here we generalize some results in the Riemannian case established in \cite{[RRS]} to the pseudo-Hermitian case.
\end{enumerate}
\end{Rk}

Taking $F(x)=x^v$ for $v>1$ and $\phi(t,x_1,x_2)=f(t)$ in the above theorem gives the corollaries as follows.
\begin{Cor} \label{thm 4.1}
Let $(M^{2m+1},H,J,\theta)$ be a pseudo-Hermitian manifold with the hypotheses \eqref{3.1}, \eqref{3.2}, \eqref{3.3}, \eqref{3.4}. Assume that $u$ is a $C^2$ function on $M$ which satisfies
\begin{align}
\triangle_b u\geq f(u)\label{4.1}
\end{align}
where $f:\mathbb{R}\rightarrow\mathbb{R}$ is a continuous function.
\begin{enumerate}[(1)]
\item If 
\begin{align}
\liminf_{x\rightarrow +\infty} \frac{f(x)}{x^v}>0\label{4.1'}
\end{align}
for some $v>1$,
then $u$ is bounded above on $M$ and $f(\sup_M u)\leq 0$.
\item If 
\begin{align}
\liminf_{x\rightarrow +\infty} \frac{f(x)}{x^v}\leq0\label{4.1''}
\end{align}
for any $v>1$,
then either $\sup_Mu=+\infty$ or $u$ is bounded above on $M$ and $f(\sup_M u)\leq 0$.
\end{enumerate}
\end{Cor}

\begin{Cor}
Let $(M^{2m+1},H,J,\theta)$ be a pseudo-Hermitian manifold with the hypotheses \eqref{3.1}, \eqref{3.2}, \eqref{3.3}, \eqref{3.4}. Assume that $u$ is a $C^2$ function on $M$ and satisfies
\begin{align}
\triangle_b u \geq c_0 u^d
\end{align}
where $d$ and $c_0$ are two constants with $d\geq 0$ and $c_0>0$. Then
\begin{enumerate}[(1)]
\item If $u\geq 0$ and $d>1$, then $u$ vanishes identically.
\item If $0\leq u\leq \sup_M u<+\infty$ and $0<d\leq 1$, then $u$ vanishes identically.
\item If $u$ is bounded above and $d=0$, then there is no such a function $u$ on $M$.
\end{enumerate}
\end{Cor}

In addition, we can also apply these generalized maximal principles to study CR conformal deformations. Before presenting the theorem, we recall that for a pseudo-Hermitian manifold $(M^{2m+1},H,J,\theta)$, if $\hat{\theta}=u^{\frac{2}{m}}\theta$ for some positive function $u\in C^\infty(M)$, then $\hat{\theta}$ is said to be conformal to $\theta$. Clearly, $\hat{\theta}$ is a pseudo-Hermitian structure of $M$ too. 

\begin{Thm}
Let $(M^{2m+1},H,J,\theta)$ be a pseudo-Hermitian manifold with the hypotheses \eqref{3.1}, \eqref{3.2}, \eqref{3.3}, \eqref{3.4} and the pseudo-Hermitian scalar curvature $\rho>0$. Let $\hat{\rho}$ be a nonpositive smooth function on $M$ satisfying 
\begin{align}
\sup_{M-K} \hat{\rho}<0
\end{align}
for some compact subset $K$ in $M$. Then the pseudo-Hermitian structure $\theta$ cannot be conformally deformed to another pseudo-Hermitian structure $\hat{\theta}$ so that the pseudo-Hermitian scalar curvature of $(M^{2m+1},H,J,\hat{\theta})$ is $\hat{\rho}$.
\end{Thm}
\proof
We will prove by contradiction. Suppose that there is a positive function $u\in C^\infty(M)$ such that $\hat{\theta}=u^{\frac{2}{m}}\theta$ and the pseudo-Hermitian scalar curvature of $(M^{2m+1},H,J,\hat{\theta})$ is $\hat{\rho}$. Then according to \cite{[HK]}, $u$ satisfies 
\begin{align}
-\frac{2m+2}{m}\triangle_b u+\rho u=\hat{\rho} u^{1+\frac{2}{m}}.
\end{align}
It is obvious that $u$ is not a constant. Since $\rho>0$, $\hat{\rho}\leq 0$ and $u>0$, we have
\begin{align}
\frac{2m+2}{m}\triangle_b u&\geq \rho u, \label{4.12}\\
\frac{2m+2}{m}\triangle_b u &\geq -\hat{\rho} u^{1+\frac{2}{m}}.\label{4.13}
\end{align}
Applying Bony's maximum principle \cite{[Bon]} to \eqref{4.12}, we see that $u$ can't attain the maximum on $M$. Let us consider the following function on $M$
\begin{align}
\phi=\frac{1}{(u+1)^a},
\end{align}
where $0<a<\frac{1}{m}$.
Since $u$ is positive, $\phi$ is also positive on $M$. Hence we apply Theorem \ref{theorem3.1} to $-\phi$, namely, for arbitrary $k\in \mathbb{N}_+$, there exists $x_k\in M$ such that
\begin{gather}
\lim_{k\rightarrow +\infty}\phi(x_k)=\inf_M{\phi}, \\
|\nabla^H \phi(x_k)|<\frac{1}{k},\label{4.2}\\
\triangle_b \phi(x_k)>-\frac{1}{k}.\label{4.3}
\end{gather}
By a direct computation, we have
\begin{gather}
\nabla^H \phi=-a\frac{\nabla^H u}{(u+1)^{a+1}},\label{4.4}\\
\triangle_b \phi=-a\frac{\triangle_b u}{(u+1)^{a+1}}+a(a+1)\frac{|\nabla^H u|^2}{(u+1)^{a+2}}\label{4.5}.
\end{gather}
From \eqref{4.13} and \eqref{4.2}-\eqref{4.5}, it follows that, at $x_k$,
\begin{align}
-\frac{1}{k}<\triangle_b \phi \leq a\frac{\hat{\rho}u^{1+\frac{2}{m}}}{(u+1)^{a+1}}+\frac{a+1}{a}(u+1)^a\frac{1}{k^2},
\end{align}
hence,
\begin{align}
-a\frac{ \hat{\rho}(x_k) u^{1+\frac{2}{m}}(x_k)}{(u(x_k)+1)^{2a+1}}<\frac{1}{k(u(x_k)+1)^a}+\frac{a+1}{ak^2},
\end{align}
where $x_k$ satisfies $\lim_{k\rightarrow \infty} u(x_k)=\sup_M u$. Since $u$ can't attain the maximum on $M$, $x_k\in M-K$ for $k$ large enough. Hence we have, for $k$ large enough,
\begin{align}
-a\frac{ \left(\sup_{M-K} \hat{\rho}\right) u^{1+\frac{2}{m}}(x_k)}{(u(x_k)+1)^{2a+1}}<\frac{1}{k(u(x_k)+1)^a}+\frac{a+1}{ak^2}.\label{5.38}
\end{align}
If $\sup_M u=\infty$, let $k\rightarrow +\infty$ in \eqref{5.38}, we get a contradiction. Consequently, $\sup_M u<\infty$. Let $k\rightarrow +\infty$ in \eqref{5.38} again, we obtain
\begin{align}
\left(\sup_{M-K} \hat{\rho}\right)\times \left(\sup_M u\right)^{1+\frac{2}{m}}\geq 0
\end{align}
which is a contradiction too, since $\sup_{M-K} \hat{\rho}<0$ and $u>0$.
\qed

\begin{Rk}
By the discussion in section \ref{section3}, we can see that for a pseudo-Hermitian manifold with pseudo-Hermitian Ricci curvature or horizontal sectional curvature $K_{\sec }^{H}$ bounded below, and $\|A\|_{C^1}$ bounded above, all results in this section are valid.
\end{Rk}

\bigskip

Yuxin Dong \ \ 

School of Mathematical Sciences

Fudan University

Shanghai, 200433, P. R. China

yxdong@fudan.edu.cn \ \ \

\bigskip

Weike Yu

School of Mathematical Sciences

Fudan University

Shanghai, 200433, P. R. China

wkyu2018@outlook.com

\end{document}